\documentclass[10pt]{amsart}
\usepackage{amsmath}
\usepackage{amsthm}
\usepackage{amssymb}
\usepackage{graphics}
\usepackage{epsfig}
\usepackage{amscd}
\usepackage{amsfonts}
\usepackage{amsbsy}
\usepackage{epsfig,afterpage}
\usepackage{psfrag}
\usepackage{epstopdf}

\newcommand{\R}{\ensuremath{\mathbb{R}}}

\newcommand{\N}{\ensuremath{\mathbb{N}}}

\newcommand{\T}{\theta}

\newcommand{\te}{\theta}

\newcommand{\s}{\ensuremath{\mathbb{S}}}
\def\p{\partial}
\def\e{\varepsilon}

\newtheorem {theorem} {Theorem}
\newtheorem {proposition} [theorem] {Proposition}

\begin{document}
\title[Limit cycles of discontinuous piecewise differential systems]
{Limit cycles of discontinuous piecewise quadratic and cubic
polynomial perturbations of a linear center}

\author[J. Llibre and Y. Tang]
{Jaume Llibre$^1$ and Yilei Tang$^2$}

\address{$^1$ Departament de Matem\`{a}tiques, Universitat Aut\`{o}noma de
Barcelona, 08193 Bellaterra, Barcelona, Catalonia, Spain}
\email{jllibre@mat.uab.cat}

\address{$^2$ School of Mathematical Sciences,  Shanghai
Jiao Tong University, Shanghai, 200240, P. R. China}
\email{mathtyl@sjtu.edu.cn}

\subjclass[2010]{34C29, 34C25, 47H11}

\keywords{Periodic solution, limit cycle, discontinuous piecewise
differential system, averaging theory}

\date{}
\dedicatory{}

\maketitle

\begin{abstract}
We apply the averaging theory of high order for computing the limit
cycles of discontinuous piecewise quadratic and cubic polynomial
perturbations of a linear center. These discontinuous piecewise
differential systems are formed by two either quadratic, or cubic
polynomial differential systems separated by a straight line.

We compute the maximum number of limit cycles of these discontinuous
piecewise polynomial perturbations of the linear center, which can
be obtained by using the averaging theory of order $n$ for
$n=1,2,3,4,5$. Of course these limit cycles bifurcate from the
periodic orbits of the linear center. As it was expected, using the
averaging theory of the same order, the results show that the
discontinuous quadratic and cubic polynomial perturbations of the
linear center have more limit cycles than the ones found for
continuous and discontinuous linear perturbations.

Moreover we provide sufficient and necessary conditions for the existence of a
center or a focus at infinity if  the discontinuous piecewise perturbations of the
linear center are general quadratic polynomials or cubic quasi--homogenous polynomials.

\end{abstract}

\section{Introduction and statement of the main Results}
\label{S-1}

The interest on the dynamics of piecewise linear differential
systems essentially started with the book of Andronov et al
\cite{AVK}, whose Russian version appeared around the 1930's. Due to
the rich dynamics of the piecewise linear differential systems, and
their applications in mechanics, electronics, economy, neuroscience,
..., these systems have been studied by researchers from many
different fields, see for instance the books of Bernardo et al
\cite{BBCK} and of Simpson \cite{Si}, the survey of Makarenkov and
Lamb \cite{ML}, and the references mentioned in all these works.

For the planar {\it continuous} piecewise linear differential
systems with two zones separated by a straight line, Lum and Chua
\cite{LC1, LC2} in 1991 conjectured that such differential systems
have at most one limit cycle. In 1998 Freire, Ponce, Rodrigo and
Torres \cite{FPRT} proved this conjecture.

While for the planar {\it discontinuous} piecewise linear
differential systems with two zones separated by a straight line Han
and Zhang \cite{HZ} obtained differential systems having two limit
cycles and conjectured that the maximum number of limit cycles of
such class of differential systems is two. Huan and Yang \cite{HY}
provided a numerical example of one of those differential system
having three limit cycles. Inspired in this numerical example Llibre
and Ponce \cite{LP} gave a proof of the existence of such three
limit cycles in the class of these differential systems. Later on
other authors also provide other discontinuous piecewise linear
differential systems with two zones separated by a straight line
also exhibiting three limit cycles, see \cite{BdM, Buz, Li2014}.
More discussion about limit cycles of discontinuous piecewise
differential systems can see references \cite{WZ2016, ZKB2006}.

Recently the averaging theory has been developed for studying the
periodic solutions of the discontinuous piecewise differential
systems. Thus Llibre, Mereu, Novaes and Teixeira \cite{LMN, LNT2015}
extended the averaging theory up to order $1$ and $2$ for studying
the periodic solutions of some discontinuous piecewise differential
systems using techniques of regularization. Later on Itikawa, Llibre
and Novaes \cite{ILN} improved the averaging theory at any order for
analyzing the periodic solutions of discontinuous piecewise
differential systems.

We consider planar discontinuous piecewise
differential systems having the line of discontinuity at $y=0$ of
the form
\begin{equation}\label{FG}
\begin{array}{l}
\dot{x}=F^{\pm}(x, y, \e),\\
\dot{y}=G^{\pm}(x, y, \e),
\end{array}
\end{equation}
where
\begin{equation*}
\begin{array}{lll}
\dot{x}=F^+(x, y, \e) &=y +\displaystyle\sum_{j=1}^{n} ~\e^j
(a_{j0}+a_{j1} x+a_{j2} y+a_{j3} x^2+a_{j4} x y+a_{j5} y^2 \\
& \qquad \qquad  +a_{j6} x^3+a_{j7} x^2 y+a_{j8} xy^2+a_{j9} y^3 ), \\
\dot{y}=G^+(x, y, \e) &=-x +\displaystyle\sum_{j=1}^{n} ~\e^j
(b_{j0}+b_{j1} x+b_{j2} y+b_{j3} x^2+b_{j4} x y+b_{j5} y^2 \\
& \qquad \qquad  +b_{j6} x^3+b_{j7} x^2 y+b_{j8} xy^2+b_{j9} y^3 ),
\end{array}
\end{equation*}
if $y\ge0$, and
\begin{equation*}
\begin{array}{lll}
\dot{x}=F^-(x, y, \e) &=y +\displaystyle\sum_{j=1}^{n} ~\e^j
(A_{j0}+A_{j1} x+A_{j2} y+A_{j3} x^2+A_{j4} x y+A_{j5} y^2 \\
& \qquad \qquad  +A_{j6} x^3+A_{j7} x^2 y+A_{j8} xy^2+A_{j9} y^3 ), \\
\dot{y}=G^-(x, y, \e) &=-x +\displaystyle\sum_{j=1}^{n} ~\e^j
(B_{j0}+B_{j1} x+B_{j2} y+B_{j3} x^2+B_{j4} x y+B_{j5} y^2 \\
& \qquad \qquad  +B_{j6} x^3+B_{j7} x^2 y+B_{j8} xy^2+B_{j9} y^3 ),
\end{array}
\end{equation*}
if $y\le0$, and where $n\in \N$, all parameters $a_{ji}, b_{ji},
A_{ji}, B_{ji}, \e \in \R$, and the perturbation parameter $|\e|$ is
small enough. Here $\mathbb N$ is the set of positive integers and
$\R$ is the set of real numbers. Notice that system \eqref{FG} is a
discontinuous piecewise differential system with the discontinuity
straight line $y=0$. As usual the dot denotes derivative with
respect to an independent real variable $t$.

In this paper we study the limit cycles of the discontinuous
piecewise quadratic (i.e. when all the cubic monomials in \eqref{FG}
are zero) and cubic polynomial differential system \eqref{FG},
which bifurcate from the periodic orbits of the linear center $\dot
x= y$, $\dot y= -x$. A classical problem for smooth differential
systems is the weak 16th Hilbert problem, which essentially asks for
the maximal number of limit cycles that bifurcate from the periodic
orbits of a center when this is perturbed inside a class of
polynomial differential systems with a fixed degree, see for more
details \cite{A1, A2, Kh, Va}. Here we are extending this problem to
the non--smooth differential system \eqref{FG}.

We denote by $L_2(n)$ and $L_3(n)$ the {\it maximum number of limit
cycles} of the discontinuous piecewise polynomial differential
system \eqref{FG} with degree $2$ and $3$ respectively {\it which
can be obtained using the averaging theory of order $n$} described
in section \ref{s2}. Then we have the following results.

\begin{theorem}\label{th1}
For $n=1,2,3,4, 5$ we have that $L_2(n)= 2, 3, 5, 6, 8$, and
$L_3(n)=3, 5, 8, 11, 13$, respectively.
\end{theorem}

Iliev in \cite{Iliev} studied the maximum number of limit cycles
$L_I(n)$ coming from the perturbation of the linear center
$\dot{x}=y, ~\dot{y}=-x$ when this center is perturbed inside the
class of all polynomial differential systems of degree $n$. Buzzi,
Pessoa and Torregrosa in \cite{Buz} found the maximum number of
limit cycles $L_1(n)$, that bifurcate from the periodic orbits of
the linear center $\dot{x}=y, ~\dot{y}=-x$ when this center is
perturbed inside the class of all discontinuous piecewise linear
differential systems separated by a straight line. Their results for
$n=1,2,3,4,5$ together with the results of Theorem \ref{th1} are
given in Table \ref{globaltable}.

\begin{table}[h]
\begin{tabular}{|c|c|c|c|c|}
  \hline

 Order $n$ & $L_1(n)$ & $L_2(n)$ & $L_3(n)$ & $L_I(n)$
 \\ \hline
   1 & 1 & 2 & 3 & 0
 \\ \hline
   2 & 1 & 3 & 5 & 1
\\ \hline
   3 & 2 & 5 & 8 & 1
 \\ \hline
   4 & 3 & 6 &11 & 2
  \\ \hline
   5 & 3 & 8 &13 & 2
  \\
  \hline
\end{tabular}
\bigskip \caption{Maximum number of limit cycles bifurcating from the
periodic orbits of the linear center using averaging theory of order
$n$.} \label{globaltable}
\end{table}

If  there exists a neighborhood of the infinity in the Poincar\'e
disc \cite[Chapter 5]{DLA2006} filled of periodic orbits, then we
say that system \eqref{FG} has a {\it center at infinity}. If there
exists a neighborhood of the infinity in the Poincar\'e disc where
all the orbits spiral going to or coming from the infinity, then we
say that system \eqref{FG} has a {\it focus at infinity}. We shall
investigate the problem of the existence of a center or a focus at
infinity under small perturbations, but before we need some
definitions.

A planar polynomial differential system
\begin{eqnarray}\label{quasi5-general}
\dot{x}= P(x, y), \qquad\qquad  \dot{y} = Q(x, y),
\end{eqnarray}
where $P(x, y)$ and $Q(x, y)$ are non--zero polynomials, is {\it
quasi--homogeneous} if there exist $s_1,s_2,d\in \mathbb N$ such
that for all positive number $\alpha$ they satisfy
\begin{equation*}  \label{alphaPQ}
P(\alpha^{s_1}x, \alpha^{s_2}y) =
\alpha^{s_1+d-1} P(x, y), \ \  \ \  Q(\alpha^{s_1}x, \alpha^{s_2}y)=
\alpha^{s_2+d-1}  Q(x, y),
\end{equation*}
Then as usual $(s_1,s_2)$ are the {\it weight exponents}, $d$ is the
{\it weight degree} with respect to the weight exponents, and $w
=(s_1, s_2, d)$ is the {\it weight vector} of the quasi--homogeneous
polynomial differential system (\ref{quasi5-general}).

By Proposition 19 of Gin\'e, Grau and Llibre \cite{Garc2013}, an
irreducible quasi-homogeneous but non-homogeneous cubic ordinary
polynomial differential system  can be written in one of the
following forms:
\begin{align*}
(I) & ~\dot{x} =y(a_1x+b_1y^2), ~\dot{y} =c_1x+d_1y^2, \mbox{\it with} ~  b_1c_1 \ne0.
\\
(II) & ~\dot{x}=a_2x^2+b_2y^3,~\dot{y}=c_2xy, ~\mbox{\it with}~ a_2b_2c_2\ne 0.
\\
(III) & ~\dot{x}=a_3y^3,~\dot{y}=b_3 x^2~\mbox{\it with}~ a_3b_3\ne 0.
\\
(IV) & ~\dot{x}=x(a_4x+b_4y^2),~\dot{y}=y(c_4x+d_4y^2),~\mbox{\it with}
~  a_4d_4\ne 0.
\\
(V) & ~\dot{x}=a_5xy^2,~\dot{y}=b_5 x^2+c_5y^3, ~\mbox{\it with}~ a_5b_5c_5\ne 0.
\\
(VI) & ~\dot{x}=a_6xy^2,~\dot{y}=b_6x+c_6y^3, ~\mbox{\it with } ~a_6b_6c_6\ne 0.
\\
(VII) & ~\dot{x}=a_7x+b_7y^3,~\dot{y}=c_7y, ~\mbox{\it with} ~a_7b_7c_7\ne 0.
\end{align*}
Perturbing the linear center by discontinuous cubic quasi-homogenous
but non-homogeneous polynomials, we obtain the following $7$ systems:
 \begin{equation}
 \label{I-piece}
\begin{array}{l}
\dot{x} =y + ~\e    y(a_1x+b_1y^2), ~~ \dot{y} =-x + ~\e
(c_1x+d_1y^2) ~\mbox{\it if} ~y\ge 0,
\\
\dot{x} =y + ~\e    y(A_1x+B_1y^2), ~~ \dot{y} =-x + ~\e
(C_1x+D_1y^2)  ~\mbox{\it if} ~y\le 0
\end{array}
\end{equation}
where 
$b_1c_1 B_1C_1 \ne0$,
 \begin{equation}
 \label{II-piece}
\begin{array}{l}
\dot{x} =y + ~\e     a_2x^2+b_2y^3, ~~ \dot{y} =-x + ~\e     c_2xy
~\mbox{\it if} ~y\ge 0,
\\
\dot{x} =y + ~\e     A_2x^2+B_2y^3, ~~ \dot{y} =-x + ~\e     C_2xy
~\mbox{\it if} ~y\le 0,
\end{array}
\end{equation}
where $a_2b_2c_2A_2B_2C_2\ne 0$,
 \begin{equation}
 \label{III-piece}
\begin{array}{l}
\dot{x} =y + ~\e      a_3y^3, ~~ \dot{y} =-x + ~\e     b_3 x^2
~\mbox{\it if} ~y\ge 0,
\\
\dot{x} =y + ~\e      A_3y^3, ~~ \dot{y} =-x + ~\e     B_3 x^2
~\mbox{\it if} ~y\le 0,
\end{array}
\end{equation}
where $a_3b_3A_3B_3\ne 0$,
 \begin{equation}
 \label{IV-piece}
\begin{array}{l}
\dot{x} =y + ~\e      x(a_4x+b_4y^2), ~~ \dot{y} =-x + ~\e
y(c_4x+d_4y^2) ~\mbox{\it if} ~y\ge 0,
\\
\dot{x} =y + ~\e      x(A_4x+B_4y^2), ~~ \dot{y} =-x + ~\e
y(C_4x+D_4y^2)  ~\mbox{\it if} ~y\le 0,
\end{array}
\end{equation}
where 
$a_4d_4A_4D_4\ne 0$,
 \begin{equation}
 \label{V-piece}
\begin{array}{l}
\dot{x} =y + ~\e      a_5xy^2, ~~ \dot{y} =-x + ~\e      b_5
x^2+c_5y^3 ~\mbox{\it if} ~y\ge 0,
\\
\dot{x} =y + ~\e      A_5xy^2, ~~ \dot{y} =-x + ~\e     B_5
x^2+C_5y^3  ~\mbox{\it if} ~y\le 0,
\end{array}
\end{equation}
where $a_5b_5c_5A_5B_5C_5\ne 0$,
 \begin{equation}
 \label{VI-piece}
\begin{array}{l}
\dot{x} =y + ~\e      a_6xy^2, ~~ \dot{y} =-x + ~\e      b_6x+c_6y^3
~\mbox{\it if} ~y\ge 0,
\\
\dot{x} =y + ~\e      A_6xy^2, ~~ \dot{y} =-x + ~\e     B_6x+C_6y^3
~\mbox{\it if} ~y\le 0,
\end{array}
\end{equation}
and $a_6b_6c_6A_6B_6C_6\ne 0$, and
\begin{equation}\label{VII-piece}
\begin{array}{l}
\dot{x} =  y + ~\e  (a_7x +b_7y^3), ~~ \dot{y} =      -x + ~\e  c_7y
~\mbox{\it if} ~y\ge 0,
\\
\dot{x} =      y + ~\e  (A_7x +B_7y^3), ~~ \dot{y} =       -x + ~\e
C_7,  ~\mbox{\it if} ~y\le 0,
\end{array}
\end{equation}
where $a_7b_7c_7A_7B_7C_7 \ne 0$.

A center is called a {\it global center } when the periodic orbits
surrounding the center filled the whole plain except the center
itself.

\begin{theorem}\label{th2}
Assume $n=1$  in system \eqref{FG}.
\begin{itemize}
\item[(i)] System \eqref{FG}  has neither centers nor foci at
infinity if the discontinuous polynomial perturbations are of degree
$2$ (i.e. if $a_{ji}=b_{ji}= A_{ji}= B_{ji}=0$ for $i=6,...,9$).

\smallskip

\item[(ii)] The unique systems from  \eqref{I-piece} to
\eqref{VII-piece} which can have  a center or a focus at infinity
are the systems \eqref{I-piece} or  \eqref{VII-piece}.

\smallskip

\item[(iii)] The infinity of system  \eqref{VII-piece} is a focus.
System  \eqref{I-piece} has a focus or a center at infinity if
$-b_1\e<0$ and $-B_1\e<0$,  and it has a  center at infinity if
$-b_1\e<0$, $-B_1\e<0$, $a_1=-2 d_1$  and $A_1=0=D_1$, which is a
global center.
\end{itemize}
\end{theorem}

\section{Averaging theory and the Descartes Theorem}\label{s2}

Using the polar coordinates $(r, \T)$ such that $x = r \cos\T$ and
$y = r \sin\te$, the differential system \eqref{FG} in these
coordinates becomes
\begin{equation}\label{PQ}
\dfrac{d r}{d \T}=\left\{\begin{array}{l}
P^+(\T,r,\e) \quad \textrm{if}\quad 0\leq\T\leq\pi,\vspace{0.2cm}\\
P^-(\T,r,\e) \quad \textrm{if}\quad -\pi\le \T \le 0,
\end{array}\right.
\end{equation}
where $P^{\pm}(\T,r,\e)=\sum_{j=1}^k\e^j P_j^{\pm}(\T,r)+\e^{k+1}
Q^{\pm}(\T,r,\e)$ with $k\in \N$, $\theta \in \mathbb{S}^1$ and
$r\in\R_+$, the functions $P_j^{\pm}:\s^1\times \R_+ \rightarrow\R$ for
$j=1, 2, \ldots,k,$ and $Q^{\pm}:\s^1\times \R_+
\times(-\e_0,\e_0)\rightarrow\R$ are analytic. Here $\e_0>0$ and
$\R_+= [0,\infty)$.

The \emph{averaged function} $f_j:\R_+\rightarrow\R$ of order $j$
for the differential equation \eqref{PQ} is defined as

\begin{equation}
\label{f}
f_j(r )=\dfrac{y^+_j(\pi,r )-y^-_j(-\pi,r )}{j!},
\quad j=1,2,\ldots,k,
\end{equation}
where $y_j^{\pm}$ for $j=1,2,3,4,5$ are
\allowdisplaybreaks
\begin{align}
y^{\pm}_1(\T,r )=& \displaystyle \int_0^{\T}P^{\pm}_1(\phi,r )d\phi,\vspace{0.1cm}\nonumber
\\
y^{\pm}_2(\T,r )=&\displaystyle \int_0^{\T}\Big(2P^{\pm}_2(\phi,r )+2\p
P^{\pm}_1(\phi,r ) y^{\pm}_1(\phi,r )\Big)d\phi, \vspace{0.1cm} \nonumber
\\
y^{\pm}_3(\T,r )=&\displaystyle \int_0^{\T}\Big( 6P^{\pm}_3(\phi,r )+6\p
P^{\pm}_{2}(\phi,r ) y^{\pm}_1(\phi,r )\nonumber\\
&+3\p^2 P^{\pm}_1(\phi,r ) y^{\pm}_1(\phi,r )^{2}+3\p P^{\pm}_1
(\phi,r )\,y^{\pm}_2(\phi,r )\Big)d\phi,\vspace{0.1cm}\nonumber
\\
y^{\pm}_4(\T,r )=&\displaystyle \int_0^{\T}\Big(24P^{\pm}_4(\phi,r )+24
\p P^{\pm}_{3}(\phi,r ) y^{\pm}_1(\phi,r )  \label{yy} \\
&+12\p^2 P^{\pm}_{2}(\phi,r ) y^{\pm}_1(\phi,r )^{2}+
12\p P^{\pm}_{2}(\phi,r ) y^{\pm}_2(\phi,r )\nonumber\\
&+12\p^2 P^{\pm}_{1}(\phi,r ) y^{\pm}_1(\phi,r )
y^{\pm}_2(\phi,r )\nonumber\\
&+4\p^3 P^{\pm}_{1}(\phi,r ) y^{\pm}_1(\phi,r )^{3}+ 4\p
P^{\pm}_{1}(\phi,r ) y^{\pm}_3(\phi,r )\Big)d\phi,\vspace{0.1cm}
\nonumber\\
y^{\pm}_5(\T,r )=&\displaystyle
\int_0^{\T}\Big(120P^{\pm}_5(\phi,r )+120\p
P^{\pm}_{4}(\phi,r ) y^{\pm}_1(\phi,r )\nonumber\\
&+60\p^2 P^{\pm}_{3}(\phi,r )
y^{\pm}_1(\phi,r )^{2}+60\p P^{\pm}_{3}(\phi,r )y^{\pm}_2(\phi,r )\nonumber\\
&+60\p^2P^{\pm}_{2}(\phi,r )y^{\pm}_1(\phi,r )  y^{\pm}_2(\phi,r )
+20\p^3 P^{\pm}_{2}(\phi,r )y^{\pm}_1(\phi,r )^{3}\nonumber\\
&+20\p P^{\pm}_{2}(\phi,r )y^{\pm}_3(\phi,r )
+20\p^2P^{\pm}_{1}(\phi,r )y^{\pm}_1(\phi,r )  y^{\pm}_3(\phi,r )\nonumber\\
&+15\p^2 P^{\pm}_{1}(\phi,r )y^{\pm}_2(\phi,r )^{2}+ 30\p^3
P^{\pm}_{1}(\phi,r ) y^{\pm}_1(\phi,r )^{2}
y^{\pm}_2(\phi,r )\nonumber\\
&+5\p^4 P^{\pm}_{1}(\phi,r ) y^{\pm}_1(\phi,r )^{4}+5\p
P^{\pm}_{1}(\phi,r ) y^{\pm}_4(\phi,r )\Big)d\phi.\nonumber
\end{align}
From \cite{ILN} we have the following result.

\begin{theorem}\label{th3}
Suppose that $j$ is the first integer such that the averaged
function $f_i= 0$ for $i=1,2,\ldots, j-1$ and $f_j\ne 0$.  If there
is $r^*\in \R_+$ such that $f_j(r^*)=0$ and $f'_j(r^*)\ne 0,$ then
for $|\e|\ne 0$ small enough there is a $2\pi$--periodic solution
$r(\T,\e)$ of \eqref{PQ} such that $r(0,\e)\to r^*$ when $\e\to 0$.
\end{theorem}

Note that the simple positive zeros of the averaged function $f_j$
provides limit cycles of the differential equation \eqref{PQ}.

We shall use the following version of the Descartes Theorem as it is
proved in \cite{Berezin}.

\begin{theorem}[Descartes theorem]\label{Dth}
Consider the real polynomial $p(x) = a_{i_1} x^{i_1} + a_{i_2}
x^{i_2} + \ldots +a_{i_r} x^{i_r}$ with $0 = i_1 < i_2 < \ldots <
i_r$. If $a_{i_j}a_{i_{j+1}} < 0,$ we say that we have a variation
of sign. If the number of variations of signs is $m,$ then the
polynomial $p(x)$ has at most $m$ positive real roots. Furthermore,
always we can choose the coefficients of the polynomial $p(x)$ in
such a way that $p(x)$ has exactly $r-1$ positive real roots.
\end{theorem}

\bigskip

\section{Proof of Theorem \ref{th1}}\label{proofth1}

We write the discontinuous piecewise {\it cubic} polynomial
differential system \eqref{FG} in polar coordinates, obtaining a
differential system $(\dot r, \dot \theta)$. After taking $\T$ as
the new independent variable we get a differential equation
$dr/d\T$, and doing Taylor series expansion of $dr/d\T$ with respect
to the variable $\e$ at $\e=0$ we obtain the differential equation
\eqref{PQ} associated to system \eqref{FG}.

Since system \eqref{FG} is a polynomial differential system, the
functions $P_j^{\pm}(\theta,r)$ and $Q_j^{\pm}(\T,r,\e)$ are analytic.
Moreover the differential equation $dr/d\T$ in the form \eqref{PQ}
is $2\pi-$periodic because the variable $\T$ appears through the
sinus and cosinus functions. In order to apply Theorem \ref{th3} to
our differential equation $dr/d\T$ it suffices to take an open
interval $\mathcal{D}=\{r:0<r<r_0\} \subset \R_+$, where the
unperturbed system can have periodic orbits $r(\T)$ such that
$r(0)=r$ with $0<r<r_0.$ Here we only give the explicit expressions
of \allowdisplaybreaks
\begin{align*}
P_1^+(a_{ij}, b_{ij}, \theta, r)=&
\frac{1}{8}\bigg(r^3\Big(4b_{19}\cos(2 \theta)-2b_{18} \sin(2
\theta)+a_{19}\sin(4 \theta)- a_{17}\sin(4 \theta)
\\
&~ -4a_{16} \cos(2 \theta)-2 a_{17} \sin(2 \theta)+b_{18}\sin(4
\theta)+a_{18}\cos(4 \theta)
\\
&~-2a_{19}\sin(2 \theta)+b_{17}\cos(4 \theta)-a_{16}\cos(4
\theta)-2b_{16}\sin(2\theta)
\\
&~-b_{16}\sin(4 \theta)-b_{19}\cos(4 \theta)-a_{18}-3
a_{16}-3b_{19}-b_{17}\Big)
\\
& ~ +r^2 \Big( 2b_{15}\sin(3 \theta)-2 a_{14} \sin(3 \theta)-2
b_{14}\cos\T-2 a_{14} \sin\T
\\
&~-2b_{13}\sin(3\theta)-6 a_{13}\cos\T +2b_{14}\cos(3
\theta)-2a_{13}\cos(3 \theta)
\\
&~+2a_{15}\cos(3 \theta)-2 b_{13}\sin\T-2
a_{15}\cos\T-6b_{15}\sin\T\Big)
\\
& ~+r\Big(-4 a_{12}\sin(2 \theta)-4a_{11}\cos(2
\theta)-4a_{11}+4b_{12}\cos(2 \theta)
\\
&~-4 b_{11} \sin(2 \theta)-4b_{12} \Big) -b_{10}\sin\T -a_{10}\cos\T
\bigg),
\\
P_2^+(a_{ij}, b_{ij}, \theta, r)=&
-a_{20}\cos\theta-b_{20}\sin\theta -\Big( b_{22}- (b_{22}-a_{21}) \cos^2\theta +(a_{22}+b_{21}) \sin\theta  \cos\theta \Big) r
\\
&~
+\Big( -b_{25} \sin\theta -(b_{24} +a_{25})\cos\theta  -(a_{24} +b_{23}-b_{25})\sin\theta  \cos^2\theta
\\
&~+  (a_{25}-a_{23}+b_{24} ) \cos^3\theta \Big) r^2
+\Big(-b_{29}- (a_{29} +b_{28})\sin\theta\cos\theta
\\
&~ +(2 b_{29}-a_{28}  -b_{27}) \cos^2\theta  + ( b_{28}+a_{29}-a_{27}-b_{26})\sin\theta\cos^3\theta
\\
& ~+(b_{27}-b_{29}-a_{26}+a_{28}) \cos^4\theta\Big) r^3 ~
- \Big(b_{10}\sin\theta+a_{10}\cos\theta
\\
&~
+(b_{12}+ ( b_{11}  +a_{12}) \sin\theta\cos\theta+(a_{11}-b_{12}) \cos^2\theta) r
+( (a_{15}+b_{14}) \cos\theta
\\
&~+b_{15} \sin\theta+(  b_{13}  +a_{14}- b_{15}) \sin\theta\cos^2\theta +(a_{13} -b_{14}-a_{15}) \cos^3\theta) r^2
\\
&~ +(b_{19} +(b_{18}+a_{19}) \sin\theta\cos\theta +(a_{18} -2 b_{19} +b_{17}) \cos^2\theta
\\
&~ +(b_{19}- a_{18}-b_{17}+a_{16}) \cos^4\theta +(a_{17}+b_{16}-a_{19}-b_{18})\sin\theta\cos^3\theta) r^3\Big)
\\
&~
 \Big(b_{10}\cos\theta-a_{10}\sin\theta+(-a_{12}+(b_{12}-a_{11})\sin\theta\cos\theta +(a_{12}+b_{11})\cos^2\theta ) r
 \\
 &~ +( -a_{15} \sin\theta+(b_{15}-a_{14}) \cos\theta +(a_{15}+b_{14}-a_{13})\sin\theta \cos^2\theta
 \\
 &~ +(b_{13}- b_{15} +a_{14})\cos^3\theta) r^2
 ~ +(-a_{19}+(b_{18} +2 a_{19}-a_{17}) \cos^2\theta
 \\
 &~+ (b_{19} -a_{18})\sin\theta\cos\theta  +(b_{17}- a_{16}-b_{19}+a_{18}) \sin\theta\cos^3\theta
 \\
 &~ +(b_{16}-b_{18}+a_{17} -a_{19}) \cos^4\theta ) r^3 \Big)/r.
\end{align*}
We omit the explicit expressions of $P_k^+(a_{ij}, b_{ij},\T, r)$ for
$k=3,4,5$ because they are quite large. Moreover, we have
\begin{align*}
P_k^-(A_{ij}, B_{ij}, \theta, r)=& P_k^+(a_{ij}, b_{ij}, \theta, r),
\end{align*}
for $k=1,2,3,4,5$.

From \eqref{yy} we compute the functions $y_j^{+}(\T,r)$ and
$y_j^{-}(\T,r)$ for $j=1,\dots,5$. After we compute the averaged
functions $f_j(r)$ for $j=1,\dots, 5$ by using formulas \eqref{f}.
Thus the averaged function of first order is
\begin{equation*}
f_1(r)=\eta_{13} r^3+\eta_{12} r^2+\eta_{11} r+\eta_{10},
\end{equation*}
where
\begin{align*}
\eta_{13}=& \frac{\pi}{8}\big(- A_{18}-B_{17}-3A_{16}-3B_{19}-
b_{17} -a_{18}-3 b_{19}-3 a_{16} \big),
\\
\eta_{12}=& \frac{2}{3}\big(-
b_{13}+A_{14}+2B_{15}-2b_{15}-a_{14}+B_{13}\big),
\\
\eta_{11}=& \frac{\pi}{2}\big(- B_{12}-a_{11} - b_{12} -A_{11}\big),
\\
\eta_{10}=& -2 b_{10}+2 B_{10}.
\end{align*}
The rank of the Jacobian matrix of the function $M_1=(\eta_{13},
\eta_{12}, \eta_{11}, \eta_{10})$ with respect to the parameters
$a_{1i}, b_{1i}, A_{1i}, B_{1i}$, $i=0,1, .., 9$ is maximal, i.e. it
is $4$. Then the coefficients $\eta_{13}, \eta_{12}$, $\eta_{11}$
and $\eta_{10}$ are linearly independent in their variables. Clearly
$f_1(r)=0$ has at most three solutions in $\mathcal{D}$. Thus, by
Theorems \ref{th3} and \ref{Dth} it follows that at most $3$ limit
cycles can bifurcate from the periodic orbits of the linear
system  using
the averaging theory of first order, and from the last part of
Theorem \ref{Dth} there are systems \eqref{FG} with three limit
cycles.

Solving $\eta_{13}$ for $A_{16}$, $\eta_{12}$ for $B_{15}$,
$\eta_{11}$ for $B_{12}$ and $\eta_{10}$  for $B_{10}$, we obtain
that $f_1(r)\equiv 0$. Applying the averaging theory of order two,
we get the second averaged function
\begin{eqnarray*}
f_2(r)=\eta_{25} r^5+\eta_{24} r^4+\eta_{23} r^3+\eta_{22}
r^2+\eta_{21} r+\eta_{20},
\end{eqnarray*}
where
\allowdisplaybreaks
\begin{align*}
\eta_{25} =&  \frac{\pi}{128} \big(18  \pi a_{16} a_{18}-2 A_{16}
B_{18} -30 A_{16} B_{16}  -6 B_{17} B_{18}  -10  B_{17} B_{16}-3
A_{18}^2 \pi
\\
& ~-18 B_{17} B_{19}  \pi-18 A_{16} A_{18}  \pi- 54 A_{16} B_{19}
\pi-6 A_{18} B_{17}  \pi-18 A_{18} B_{19}  \pi
\\
& ~ +18  \pi a_{18} b_{19}+6 \pi a_{18} b_{17}-3  B_{17}^2  \pi-
27B_{19}^2 \pi-27  A_{16}^2  \pi+3  \pi a_{18}^2+27 \pi a_{16}^2
\\
& ~-6 b_{19} b_{16}+27   \pi b_{19}^2+3  \pi b_{17}^2-2 B_{17}
A_{17}  +2 B_{17} A_{19}  +6 a_{18} a_{17}  -2 a_{18} b_{16}
\\
& ~+6 a_{16} a_{19}  +2 b_{17} a_{19}  -2 b_{17} a_{17}  -10  b_{17}
b_{16}  -6 b_{17} b_{18}  +2 b_{19} a_{17}+2 a_{18} b_{18}
\\
& ~+ 54  \pi a_{16} b_{19} +18  \pi a_{16} b_{17}-10  b_{19} b_{18}
+10  a_{18} a_{19}  +30 b_{19} a_{19}  -30  a_{16} b_{16}
\\
& ~-2 a_{16} b_{18}  +10  a_{16} a_{17}  +18  \pi b_{17} b_{19}-18
A_{16} B_{17}  \pi-2 A_{18} B_{16}  +6 A_{18} A_{17}
\\
& ~+10  A_{18} A_{19}  -6 B_{19} B_{16} +2 B_{19} A_{17}  +2 A_{18}
B_{18}  +10  A_{16} A_{17}  +6 A_{16} A_{19}
\\
& ~+30 B_{19} A_{19}  -10  B_{19} B_{18}   \big),
\\
\eta_{24}=& \frac{1}{2} \big( 16 B_{14} B_{17}/45-8 a_{17}
b_{15}/15+8 a_{18} b_{14}/45+16 a_{19} b_{15}/15-4 b_{13} b_{16}/3
\\
& ~-8 A_{15} A_{16}/15+8 B_{13} B_{18}/15-16b_{15} b_{16}/15 -16
b_{14} b_{17}/45-8 b_{13} b_{18}/15
\\
& ~+4 B_{13} B_{16}/3-8 A_{15} B_{17}/45+8 A_{17} B_{15}/15+8 A_{17}
B_{13}/15+8 a_{14} a_{19}/15
\\
& ~ -4 A_{13} B_{17}/45-28 A_{13} A_{16}/15-8A_{13} A_{18}/9 -8
a_{14} b_{16}/15+4 a_{13} b_{17}/45
\\
& ~+4A_{16} B_{14}/15-32 A_{15} B_{19}/15-16 A_{19} B_{15}/15+16
B_{15} B_{16}/15
\\
& ~+8 a_{15} a_{16}/15+8 a_{13} b_{19}/5+4 a_{14} a_{17}/15-8 A_{13}
B_{19}/5+8 A_{14} B_{16}/15
\\
& ~-8 A_{18} B_{14}/45-8 A_{14} A_{19}/15+16 B_{15} B_{18}/15-32
A_{15} A_{18}/45
\\
& ~+32 a_{15} a_{18}/45+28 a_{13} a_{16}/15+8 a_{13} a_{18}/9+8
a_{15} b_{17}/45+32 a_{15} b_{19}/15
\\
& ~-4 a_{16} b_{14}/15-16 b_{15} b_{18}/15-8 a_{17} b_{13}/15+5 \pi
a_{18} b_{13}/12-4 A_{14} A_{17}/15
\\
& ~+5 \pi b_{15} b_{19}/2+5 \pi b_{13} b_{17}/12+5 \pi a_{14}
a_{16}/4+5\pi a_{18} b_{15}/6
\\
& ~+5 \pi b_{13} b_{19}/4+5 \pi b_{15} b_{17}/6+5 \pi a_{14}
a_{18}/12+5 A_{16} B_{13} \pi/4
\\
& ~+5 A_{18} B_{15} \pi/6+5 A_{14} B_{17} \pi/12+5 A_{14} A_{18}
\pi/12+5 A_{16} B_{15} \pi/2
\\
& ~+5A_{14} B_{19} \pi/4 +5 A_{14} A_{16} \pi/4+5 B_{15} B_{19}
\pi/2+5A_{18} B_{13} \pi/12
\\
& ~+5 B_{15} B_{17} \pi/6+5 B_{13} B_{19} \pi/4+5 B_{13} B_{17}
\pi/12+5 \pi a_{14} b_{17}/12
\\
& ~+5 \pi a_{14} b_{19}/4+5 \pi a_{16} b_{13}/4+5 \pi a_{16}
b_{15}/2 \big),
\\
\eta_{23} =&   \frac{1}{2} \big(32 a_{14} b_{15}/9-32 B_{15}^2/9-8
A_{14}^2/9-8 B_{13}^2/9+16 a_{14} b_{13}/9- B_{27} \pi/4
\\
& ~-32 B_{13} B_{15}/9+32 b_{13} b_{15}/9-16 A_{14} B_{13}/9-32
A_{14} B_{15}/9- A_{28} \pi/4
\\
& ~-3 B_{29} \pi/4-3 A_{26} \pi/4-b_{27} \pi/4 - a_{28} \pi/4-3
a_{26} \pi/4+3 \pi^2 a_{16} b_{12}/4
\\
& ~ -3 b_{29} \pi/4+32 b_{15}^2/9 -\pi a_{11} b_{18}/16+3\pi^2
a_{11} b_{19}/4 -11A_{11} B_{16} \pi/16
\\
& ~+ \pi^2 a_{18} b_{12}/4- A_{11} B_{17} \pi^2/4 -3A_{11} B_{19}
\pi^2/4 -3A_{16} B_{12} \pi^2/4
\\
& ~ -3 A_{11} A_{16} \pi^2/4-3 B_{12} B_{19} \pi^2/4- B_{12} B_{17}
\pi^2/4- A_{18} B_{12} \pi^2/4
\\
& ~ -A_{11} A_{18} \pi^2/4 +A_{12} B_{17} \pi /8 + A_{18} B_{11}
\pi/8 - B_{14} B_{13} \pi/4+8 a_{14}^2/9
\\
& ~- A_{11} B_{18} \pi/16-7 A_{17} B_{12} \pi/16+3 A_{18} A_{12}
\pi/8+3 A_{11} A_{19} \pi/16
\\
& ~+ \pi^2 a_{11} a_{18}/4+ \pi^2 a_{11} b_{17}/4 +3a_{16} a_{12}
\pi/8 -b_{17} b_{11} \pi/8+8b_{13}^2/9
\\
& ~+\pi^2 b_{12} b_{17}/4 +3\pi^2 b_{12} b_{19}/4 + \pi a_{12}
b_{17}/8+ \pi a_{18} b_{11}/8+A_{14} A_{13} \pi/4
\\
& ~-7 \pi a_{17} b_{12}/16-11 b_{12} b_{16} \pi/16+3 A_{16} A_{12}
\pi/8+ A_{11} A_{17} \pi/16
\\
& ~+A_{14} A_{15} \pi/4 - A_{13} B_{13} \pi/2+9 b_{19} a_{12}
\pi/8-b_{14} b_{15} \pi/4 +3 b_{12} a_{19} \pi/16
\\
& ~-11 B_{12} B_{16} \pi/16+3 B_{12} A_{19} \pi/16-B_{17} B_{11}
\pi/8+3 B_{19} B_{11} \pi/8
\\
& ~ -9 B_{12} B_{18} \pi/16+9 B_{19} A_{12} \pi/8+3 a_{18} a_{12}
\pi/8-11 a_{11} b_{16} \pi/16
\\
& ~+a_{11} a_{17} \pi/16-3a_{16} b_{11} \pi/8 +3 b_{19} b_{11}
\pi/8+a_{15} a_{14} \pi/4-B_{14} B_{15} \pi/4
\\
& ~+a_{15} b_{15} \pi/2 -b_{14} b_{13} \pi/4 -b_{13} a_{13} \pi/2
+a_{14} a_{13} \pi/4 -9 b_{12} b_{18} \pi/16
\\
& ~-3A_{16} B_{11} \pi/8+ B_{15} A_{15} \pi/2+3\pi^2 a_{11}
a_{16}/4+3a_{11} a_{19} \pi/16\big),
\\
\eta_{22}=&  \frac{1}{2} \big(4 B_{23}/3-4 a_{24}/3-4b_{23}/3
-8b_{25}/3 +4A_{24}/3 +8B_{25}/3-4b_{11} b_{13}/3
\\
& ~ -4A_{12} A_{14}/3 -16b_{12} b_{14}/9 +8a_{12} b_{15}/3 -20
a_{13} b_{12}/9-8A_{12} B_{15}/3
\\
& ~+8a_{11} a_{15}/9 -4a_{11} b_{14}/9 -8A_{15} B_{12}/9 -8a_{17}
b_{10}/3 +4a_{11} a_{13}/9
\\
& ~ -8b_{10} b_{18}/3 -4 b_{10} b_{16}+8 a_{10} b_{19}+4a_{12}
a_{14}/3 +8a_{10} a_{18}/3 +4a_{10} b_{17}/3
\\
& ~+8A_{17} B_{10}/3 +16B_{12} B_{14}/9 -4A_{10} B_{17}/3 -8 A_{10}
B_{19}+4 B_{11} B_{13}/3
\\
& ~ -8A_{10} A_{18}/3 +8 B_{10} B_{18}/3+4 B_{10} B_{16}+4A_{11}
B_{14}/9 +20A_{13} B_{12}/9
\\
& ~-8A_{11} A_{15}/9 -4A_{11} A_{13}/9 -4 A_{10} A_{16}+8a_{15}
b_{12}/9 +4 a_{10} a_{16}+\pi a_{11} a_{14}
\\
& ~+\pi a_{11} b_{13}+\pi a_{14} b_{12}+\pi b_{12} b_{13}+3a_{18}
b_{10} \pi/4 +9a_{16} b_{10} \pi/4 +3b_{10} b_{17} \pi/4
\\
& ~+2 \pi b_{12} b_{15}+9b_{10} b_{19} \pi/4 +2 \pi a_{11}
b_{15}+B_{12} B_{13} \pi+A_{14} B_{12} \pi+A_{11} B_{13} \pi
\\
& ~+A_{11} A_{14} \pi+3B_{10} B_{17} \pi/4 +3A_{18} B_{10} \pi/4
+9A_{16} B_{10} \pi/4 +9B_{10} B_{19} \pi/4
\\
& ~+2 B_{12} B_{15} \pi+2 A_{11} B_{15} \pi \big),
\\
\eta_{21}=& \frac{1}{2} \big(16 b_{10} b_{15}/3-16B_{10} B_{15}/3+8
b_{10} b_{13}/3+8a_{14} b_{10}/3 - B_{10} B_{13}8/3
\\
& ~ -8 A_{14} B_{10}/3-B_{22} \pi- A_{11}^2 \pi^2/4-B_{12}^2 \pi^2/4
+\pi^2 b_{12}^2/4 + \pi^2 a_{11}^2/4
\\
& ~-A_{21} \pi-a_{21} \pi-b_{22} \pi+2 B_{15} A_{10} \pi+a_{10}
a_{14} \pi+2 a_{10} b_{15} \pi-b_{14} b_{10} \pi
\\
& ~-2 b_{10} a_{13} \pi- b_{12} b_{11} \pi/2+ a_{11} a_{12} \pi/2+
b_{12} a_{12} \pi/2+\pi^2 a_{11} b_{12}/2
\\
& ~-a_{11} b_{11} \pi/2 +A_{11} A_{12} \pi/2 -B_{14} B_{10} \pi-2
A_{13} B_{10} \pi+B_{12} A_{12} \pi/2
\\
& ~-B_{12} B_{11} \pi/2 +A_{14} A_{10} \pi- A_{11} B_{12} \pi^2/2-
A_{11} B_{11} \pi/2 \big),
\\
\eta_{20}=& -2 b_{20}+2 B_{20}-2 A_{10} A_{11}+ b_{10} b_{12} \pi/2+
a_{11} b_{10} \pi/2+A_{11} B_{10} \pi/2
\\
& ~-2 A_{10} B_{12}+2 a_{10} a_{11}+ B_{10} B_{12} \pi/2+2 a_{10}
b_{12}+2 B_{10} B_{11}-2 b_{10} b_{11}.
\end{align*}
Because the rank of the Jacobian matrix of the function
$M_2=(\eta_{25}, \eta_{24},\eta_{23}, \eta_{22}, \eta_{21}$,
$\eta_{20})$ with respect to its variables $a_{li}, b_{li}, A_{li},
B_{li}$, $l=1,2, i=0,1, .., 9$ is maximal, i.e. it is $6$, the
functions $\eta_{25}, \eta_{24}, \eta_{23}, \eta_{22}, \eta_{21}$
and $\eta_{20}$ are linearly independent in their variables. Hence
by Theorem \ref{Dth}, the equation $f_2(r)=0$ has at most $5$ roots
in $\mathcal{D}$ and therefore at most $5$ limit cycles of system
\eqref{FG} can bifurcate from the periodic orbits of the linear
system using the averaging theory of order two, and there are
systems \eqref{FG} having $5$ limit cycles.

Solving  $\eta_{25}, \eta_{24}, \eta_{23}, \eta_{22}, \eta_{21},
\eta_{20}$ for $A_{21}, A_{24},B_{20},  B_{27}, A_{17}, a_{17}$, we
get $f_2(r)\equiv 0$, and we can use the averaging theory of order
three. Then the third averaged function is
\begin{equation*}
rf_3(r)=\eta_{38} r^8+\eta_{37} r^7+\eta_{36} r^6+\eta_{35}
r^5+\eta_{34} r^4+\eta_{33} r^3+\eta_{32} r^2+\eta_{31} r+\eta_{30}.
\end{equation*}
The functions $\eta_{3j}$ for $j=0,\dots, 8$ are linearly
independent in their variables, because the rank of the Jacobian
matrix $M_3=(\eta_{30},\dots, \eta_{38})$ with respect to its
variables is maximal, i.e. it is $9$. We do not provide their
explicit expressions, because they are very long. Therefore the
equation $f_3(r)=0$ has at most $8$ zeros in $\mathcal{D}$ and at
most $8$ limit cycles of system \eqref{FG} can bifurcate from the
periodic orbits of the linear system using the averaging theory of
order three, and again there systems \eqref{FG} with $8$ limit
cycles.

By choosing conveniently some variables to cancel the coefficients
$\eta_{3j}$ for $j=0,\dots, 8$ we do the third order averaged
function identically zero. So we can compute the fourth averaged
function $f_4(r)$. And by doing this $f_4(r)$ identically zero we
also can compute the fifth averaged function $f_5(r)$. These two
averaged functions have the form
\begin{eqnarray*}
r^2f_4(r)&=& \eta_{411} r^{11}+\eta_{410} r^{10}+\eta_{49}
r^9+\eta_{48} r^8+\eta_{47} r^7+\eta_{46} r^6
\\
&& ~+\eta_{45} r^5+\eta_{44} r^4+\eta_{43} r^3+\eta_{42}
r^2+\eta_{41} r+\eta_{40},
\\
r^2f_5(r)&=& \eta_{513} r^{13}+\eta_{512} r^{12}+\eta_{511}
r^{11}+\eta_{510} r^{10}+\eta_{59} r^9+\eta_{58} r^8+\eta_{57} r^7
\\
&& ~+\eta_{56} r^6+\eta_{55} r^5+\eta_{54} r^4+\eta_{53}
r^3+\eta_{52} r^2 +\eta_{51} r+\eta_{50}.
\end{eqnarray*}
We can prove that the coefficients $\eta_{ij}$ are linearly
independent in their variables. Their expressions are very long so
we do not give them here. As a result of these calculations it
follows that $f_4(r)=0$ (resp. $f_5(r)=0$) has at most $11$ (resp.
$13$) solutions in $\mathcal{D}$, and therefore at most $11$ (resp.
$13$) limit cycles of system \eqref{FG} can bifurcate from the
periodic orbits of the linear center, and there are systems
\eqref{FG} having $11$ (resp. $13$) limit cycles.

Now we consider the discontinuous piecewise {\it quadratic}
polynomial perturbations in system \eqref{FG}. Doing in the previous
averaged functions $f_k(r)$ for $k=1,2,3,4,5$ the coefficients
$a_{ji}=b_{ji}=A_{ji}=B_{ji}=0$ for $i=6,...,9$, we obtain the
averaged functions for the quadratic polynomial perturbations in
system \eqref{FG}. From these averaged functions we obtain the
numbers $L_2(n)$ for $n=1,2,3,4,5$ in Theorem \ref{th1}. This
completes the proof of Theorem \ref{th1}.

\bigskip

\section{Proof of Theorem \ref{th2}}

Consider system  \eqref{FG} having the  linear center $\dot{x}=y,
~\dot{y}=-x$  and being perturbed inside the class of discontinuous
piecewise quadratic polynomial differential systems
\begin{equation}\label{deg2-1}
\begin{array}{l}
\dot{x} =y + ~\e  F_1^{\pm} (x,y), \vspace{0.2cm}
\\
\dot{y} =-x + ~\e  G_1^{\pm} (x,y),
\end{array}
\end{equation}
where
\begin{align*}
F_1^{\pm} (x,y) &=  a_0^{\pm}+a_1^{\pm} x+a_2^{\pm} y+a_3^{\pm} x^2+a_4^{\pm} x y+a_5^{\pm} y^2,
\\
G_1^{\pm} (x,y) &= b_0^{\pm}+b_1^{\pm} x+b_2^{\pm} y+b_3^{\pm} x^2+b_4^{\pm} x y+b_5^{\pm} y^2,
\end{align*}
are defined in the regions $\{y\ge0\}$ and $\{y\le 0\}$, and all
parameters $a_j^{\pm}, b_j^{\pm} \in \R$ for $j=0,1,.., 5$. It is not
difficult to find that systems \eqref{deg2-1} have not a center at
infinity because first the equation
\[
x (b_3^{\pm} x^2+b_4^{\pm} x y+b_5^{\pm} y^2) -y (a_3^{\pm} x^2+a_4^{\pm} x y+a_5^{\pm}
y^2)=0
\]
has at least a real solution because it is a cubic homogeneous
polynomial, and therefore system \eqref{deg2-1} has singularities at
infinity.   Moreover, by the analysis of the local phase portraits
of the infinite singularities of quadratic systems in
\cite{Coll1987} or \cite{SV2005}, it follows that the infinity of
system \eqref{deg2-1} cannot be a center or a focus, because always
some orbits have their $\alpha$-- or $\omega$--limits at some
infinite singularity. Hence statement (i) of Theorem \ref{th2} is
proved.

Perturbing the linear center by discontinuous cubic quasi-homogenous
but non-homogeneous polynomials, we obtain
\begin{equation}\label{deg3-1}
\begin{array}{l}
\dot{x} =y + ~\e  F_2^{\pm} (x,y),
\\
\dot{y} =-x + ~\e  G_2^{\pm} (x,y),
\end{array}
\end{equation}
where $ F_2^{\pm} (x,y)$ and $G_2^{\pm} (x,y)$  belonging to one of systems
$(I)-(VII)$ in the Section \ref{S-1},  are defined in the regions
$\{y\ge 0\}$ and $\{y\le 0\}$.   Notice that system  \eqref{deg3-1}
has not  a center or a focus at infinity if $ F_2^{\pm} (x,y)$ and
$G_2^{\pm} (x,y)$ have the forms $(III)- (VI)$  because  one of
singularities at infinity of these systems is a  saddle, node,
saddle-node or a nilpotent equilibrium by Poincar\'e transformations
\begin{equation}
x=1/z, \, y=u/z, ~\mbox{ and}~   x=v/z, \, y=1/z
\label{change-uvz}
\end{equation}
together with the time variables $d\tau=dt/z^2$. For simplicity, we
only give the  compactification systems for \eqref{deg3-1} when $
F_2^{\pm} (x,y)$ and $G_2^{\pm} (x,y)$  belonging to system $(IV)$. System
\begin{equation}\label{Per-IV}
\begin{array}{l}
\dot{x} =y + ~\e   x(a_4x+b_4y^2),
\\
\dot{y} =-x + ~\e   y(c_4x+d_4y^2),
\end{array}
\end{equation}
around the equator of the Poincar\'e sphere can be written
respectively in
\begin{equation*}\label{IV-uz}
\begin{array}{l}
\dot{u} = - \e (a_4-c_4) u z-z^2-\e (b_4-d_4) u^3-u^2 z^2,
\\
\dot{z} = -z ( \e a_4 z+\e b_4 u^2+u z^2),
\end{array}
\end{equation*}
and
\begin{equation}\label{IV-vz}
\begin{array}{l}
\dot{u} =  \e (b_4-d_4) v+z^2+\e (a_4-c_4)  v^2 z+  v^2 z^2,
\\
\dot{z} = z (-\e d_4 - \e c_4 v z+v z^2),
\end{array}
\end{equation}
after changes \eqref{change-uvz}, where  $ a_4d_4\ne 0 $. Notice
that the origin of \eqref{IV-vz}, which  is located at the end of
the $y$--axis and is a singularity at infinity of system
\eqref{Per-IV} of hyperbolic type if $b_4-d_4\ne0$ or  of
semi-hyperbolic type if $b_4-d_4=0$. Then by Theorems 2.15 and 2.19
of \cite{DLA2006} this singularity  can only be  a  saddle, node or
saddle-node. Hence there exist no centers or foci at infinity of
systems  \eqref{Per-IV}.

When $ F_2^{\pm} (x,y)$ and $G_2^{\pm} (x,y)$ have the forms $(I), (II)$,
or $(VII)$, systems \eqref{deg3-1} becomes
\begin{equation}\label{Per-I}
\dot{x} =y + ~\e    y(a_1x+b_1y^2),
~~
\dot{y} =-x + ~\e    (c_1x+d_1y^2),
\end{equation}
\begin{equation}\label{Per-II}
\dot{x} =y + ~\e   ( a_2x^2+b_2y^3),
~~
\dot{y} =-x + ~\e   c_2xy,
\end{equation}
\begin{equation}\label{Per-VII}
\dot{x} =y + ~\e  (a_7x+b_7y^3),
~~
\dot{y} =-x + ~\e  c_7y,
\end{equation}
respectively, which have no  singularities at infinity at the
endpoints of the $y$-axis from simple calculations. By the first
change of  \eqref{change-uvz}, systems \eqref{Per-I}-\eqref{Per-VII}
can be transformed into
\begin{equation}\label{I-uz}
\begin{array}{l}
\dot{u} =  (c_1 \e-1) z^2-u^2z^2- \e (a_1-d_1) u^2 z-b_1 \e u^4,
\\
\dot{z} =  -u z (a_1 \e z+ \e b_1u^2+z^2),
\end{array}
\end{equation}
\begin{equation}\label{II-uz}
\begin{array}{l}
\dot{u} =   - \e (a_2-c_2) uz - z^2-u^2z^2-b_2 \e u^4,
\\
\dot{z} =  -z (\e a_2 z+\e b_2 u^3+u z^2),
\end{array}
\end{equation}
and
\begin{equation}\label{VII-uz}
\begin{array}{l}
\dot{u} =   -z^2  -\e (a_7-c_7) u z^2-\e b_7u^4-u^2z^2 ,
\\
\dot{z} =  -z (\e a_7 z^2+\e b_7u^3+u z^2),
\end{array}
\end{equation}
respectively. Note that all systems \eqref{I-uz}-\eqref{VII-uz} have
a unique singularity at the end points of  the $x$-axis, which
corresponds to the  origin denoted by $C_j=(0,0)$ for $j=I, II,
VII$. We will analyze the local phase portrait of the singularity
$C_j$.

First we consider the properties of $C_I$ for system \eqref{I-uz}.
Notice that the vector field of \eqref{I-uz} is invariant under the
change of variables $(u, z, t)\to (-u, z, -t)$. Therefore, system
\eqref{I-uz} is symmetric with respect to the $z$-axis. Thus, we
only need to consider the right half-plane $u\ge0$ for studying the
local phase portrait of $C_I$. Using the change $u_1=u^2,z_1=z$,
system \eqref{I-uz} becomes
\begin{equation}\label{I-uz2}
\begin{array}{l}
\dot{u}_1 = -2 \e b_1 u_1^2-2 \e (a_1-d_1) u_1 z+  2 (c_1 \e-1) z^2-u_1z^2,
\\
\dot{z} = -z (\e b_1u_1+\e a_1 z+z^2).
\end{array}
\end{equation}

We need to use the following notions. Consider the analytic
differential system
\begin{eqnarray} \label{qu.deg}
\left.  \begin{array}{ll}
   \dot{x}    =X_m(x,y)+\Phi_m(x,y):=X(x,y),
   \\
   [1ex]
   \dot{y}    =Y_m(x,y)+\Psi_m(x,y):=Y(x,y),
  \end{array}
  \right.
  \end{eqnarray}
where $X_m(x,y)$ and $Y_m(x,y)$ are homogeneous polynomials of
degree $m\ge 1$ such that simultaneously do not vanish,  and
$\Phi_m(x,y),\,\, \Psi_m(x,y)=o(r^m)$ when $r=\sqrt{x^2+y^2} \to 0$.
Let the origin $O$ be an isolated singularity of \eqref{qu.deg}. In
order to see when there exist orbits connecting with   $O$, by
Lemmas 1 and 3 in \cite[Chapter 2]{Sans} we only need to discuss the
orbits along exceptional directions of system \eqref{qu.deg} at $O$.

Applying the polar coordinate changes $x=r \cos\theta$ and $y=r
\sin\theta$, system \eqref{qu.deg} can be written in
\begin{eqnarray}
 \frac{1}{r}\frac{dr}{d\theta}=\frac{H_m(\theta)+o(1)}{G_m(\theta)+o(1)},
\quad {\rm as}\;\; r\to 0, \label{H/G}
\end{eqnarray}
where
\begin{eqnarray*}
  \begin{array}{ll}
   G_m(\theta)=\cos\theta Y_m(\cos\theta,\sin\theta)-
            \sin\theta X_m(\cos\theta,\sin\theta).
     \\
   H_m(\theta)=\sin\theta Y_m(\cos\theta,\sin\theta)+
            \cos\theta X_m(\cos\theta,\sin\theta).
  \end{array}
\label{GH0}
\end{eqnarray*}
Hence a necessary condition for $\theta=\theta_0$ to be an
exceptional direction is $G_m(\theta_0)=0$.

For our system \eqref{I-uz2} we calculate
\[
G_2(\theta)=\sin\theta \Big( b_1 \e \cos^2\theta+(-a_1 \e+2 \e
(a_1-d_1)) \sin\theta\cos\theta+(-2 c_1 \e+2) \sin^2\theta \Big).
\]
When $-b_1\e<0$, the equation $G_2(\theta)=0$ has only two zeros $0$
and $\pi$ if $\theta\in[0,2\pi)$. When $\theta\to 0$ equation
\eqref{H/G} has the form
\[
\frac{1}{r}\frac{dr}{d\theta}=\frac{H_2(\theta)+o(1)}{G_2(\theta)+o(1)}=
-\frac{2}{\theta}+O(1).
\]
Then $r=r_1 e^{\int_{\theta_0}^{\theta} ~-\frac{2}{\theta}+O(1)}
d\theta \to +\infty$ as $\theta\to 0$. Thus, by a similar proof of
Theorem 10.1 and Theorem 10.5 in \cite{Ye1986}, we obtain that  the
$u$-axis is the only orbit connecting with the origin of system
\eqref{I-uz2} if $-b_1\e<0$. Therefore the discontinuous piecewise
cubic polynomial differential system  \eqref{I-piece} has  a center
or a focus at infinity when $-b_1\e<0$ and $-B_1\e<0$, where $
b_1c_1  B_1C_1 \ne0$.

Second we consider the local phase portrait of $C_{II}=(0,0)$ for
system \eqref{II-uz}. By the blow-up $u=u_2 z$ along the $u$-axis
together with a time scaling $dt= dt_1/z$, system \eqref{II-uz}
becomes
\begin{equation}\label{II-uz-BU}
\begin{array}{l}
\dot{u} =   -1+ c_2 \e u,
\\
\dot{z} =   -z (a_2 \e+u z^2+ b_2\e u^3 z^2),
\end{array}
\end{equation}
where we still write $u_2$ as $u$ for simplicity. The $u$-axis
system \eqref{II-uz-BU} has only one singularity $(1/(c_2 \e), 0)$,
which is a saddle or node if $a_2c_2 \e\ne0$. Therefore there exist
neither centers nor foci at infinity of systems  \eqref{II-piece}.

Finally we consider the local phase portrait of $C_{VII}=(0,0)$ for
system \eqref{VII-uz}. By the blow-up $u=u_7 z$ along the $u$-axis
together with a time scaling $dt= dt_7/z$, system \eqref{VII-uz}
becomes
\begin{equation}\label{VII-uz-BU}
\begin{array}{l}
\dot{u} =  -1 + c_7 \e u z,
\\
\dot{z} =   -z^2 (a_7 \e +u z +b_7 \e u^3 z),
\end{array}
\end{equation}
where we still write $u_7$ as $u$ for simplicity. The $u$-axis
system \eqref{VII-uz-BU} has no singularities. Notice that the
$u$-axis is an orbit of system \eqref{VII-uz-BU}  and no other
orbits can connect with the $u$-axis. Furthermore for system
\eqref{VII-uz} we calculate $ \tilde{G}_2(\theta)=\sin^3\theta $ in
\eqref{H/G}, implying that all possible orbits connecting with
$C_{VII}$ must be along the direction of $u$-axis. Therefore we
obtain that no orbits can go to or come from the singularities at
infinity of system \eqref{Per-VII}, consequently no orbits can go to
or come from the singularities at infinity of system
\eqref{VII-piece}

This completes the proof of statement (ii) in  Theorem \ref{th2}.

\smallskip

Now we shall study the existence or not of a center or a focus at
infinity for the discontinuous piecewise cubic polynomial
differential systems \eqref{I-piece} and \eqref{VII-piece}, for this
we shall use the averaging theory for proving statement (iii) in
Theorem \ref{th2}.

We claim that the infinity  is a focus for system \eqref{VII-piece}.
In fact, from the formula \eqref{f}, we can compute the averaged
functions of order $\le4$ for system \eqref{VII-piece}. The averaged
function of first order is
\begin{equation*}
\tilde{f}_1(r)= - r \pi (a_7+c_7+A_7+C_7)/2.
\end{equation*}
Then we solve $a_7 = -c_7-A_7-C_7$ from $\tilde{f}_1(r)\equiv 0$.
Applying the averaging theory of order two,
we get the second averaged function
\begin{eqnarray*}
\tilde{f}_2(r)= 3 \pi r^3 (B_7-b_7) (A_7+C_7)/32.
\end{eqnarray*}
Solving $\tilde{f}_2(r)\equiv 0$, we obtain $b_7= B_7 $ or $A_7=
-C_7$. We use the averaging theory of order three and get the third
averaged function
\begin{equation*}
\tilde{f}_3(r)= \pi(A_7+C_7) r\Big(-\frac{35}{512} (B_7^2-b_7^2)
r^4+ \frac{3}{32} \pi(B_7-b_7) (A_7+C_7) r^2+ \frac{1}{4} (C_7+c_7)
(A_7+c_7)  \Big).
\end{equation*}
When $A_7= -C_7$  we have $\tilde{f}_3(r)\equiv 0$, and when $b_7=
B_7 $ we have $\tilde{f}_3= \pi (C_7+c_7) (A_7+c_7) (A_7+C_7) r /4$.
Hence we get $(C_7+c_7) (A_7+c_7) (A_7+C_7)=0$ from
$\tilde{f}_3(r)\equiv 0$. We compute the averaged function of order
$4$ for system \eqref{VII-piece} and we get
\begin{equation*}
\tilde{f}_{41}(r)= B_7 \pi r^3 \Big( - \frac{3}{128} (2 A_7+3 C_7)
(-C_7+A_7)^2 -  \frac{33}{1024}  B_7^2   (2 A_7+5 C_7) r^4  \Big)
\end{equation*}
if $b_7= B_7$ and $ A_7+c_7=0$, or
\begin{equation*}
\tilde{f}_{42}(r)=  \frac{3}{1024} B_7 \pi r^3\Big( - (16 A_7^3-8
A_7^2 C_7-32 A_7 C_7^2+24 C_7^3)- (22 A_7 B_7^2+55 B_7^2 C_7) r^4
\Big)
\end{equation*}
if $b_7= B_7$ and $ C_7+c_7=0$,  or
\begin{equation*}
\tilde{f}_{43}(r)= B_7 C_7\pi r^3 \Big( -\frac{99}{1024} B_7^2
r^4-\frac{3}{32} C_7^2 \Big)
\end{equation*}
if $A_7= -C_7$. Therefore, the functions $\tilde{f}_{41}(r)$,
$\tilde{f}_{42}(r)$ and $\tilde{f}_{43}(r)$ cannot be identically
zero, otherwise we have a contradiction with the fact that $B_7
C_7\ne 0$. So there are no periodic orbits, and we can obtain some
isolated spiral orbit as close as we want to infinity. Therefore the
infinity is a focus for system \eqref{VII-piece} and the claim is
proved.

\smallskip

\begin{figure}[htb]
\centering
\includegraphics[scale=.60]{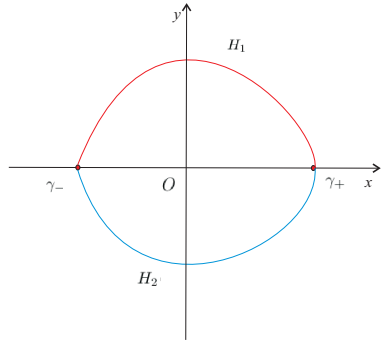}
\caption{  Existence of closed orbits for system \eqref{I-piece2}.  }
\label{Center1}
\end{figure}

For system \eqref{I-piece} we can compute the averaged functions of
order $\le5$. Here we omit the tedious calculations and only show
the results. Doing all the averaged functions of order less than $5$
identically zero, we obtain $a_1=-2 d_1$ and $A_1=0=D_1$. Then
system \eqref{I-piece} becomes
\begin{equation}\label{I-piece2}
\begin{array}{lll}
&\dot{x} =y + ~\e    y(-2 d_1 x+b_1y^2), ~~
&\dot{y} =-x + ~\e    (c_1x+d_1y^2), ~\mbox{\it if} ~y\ge 0,
\\
&\dot{x} =y + ~\e    B_1y^3, ~~
&\dot{y} =-x + ~\e    C_1x,  ~\mbox{\it if} ~y\le 0,
\end{array}
\end{equation}
where $-b_1\e<0$, $-B_1\e<0$ and $ b_1c_1 B_1C_1 \ne0$. System
\eqref{I-piece2} has the polynomial first integral
\[
H_1(x, y)=\frac{y^2}{2} +\frac{\e b_1}{4} y^4 + \frac{(1-\e
c_1)x^2}{2} -\e d_1 xy^2
\]
if $y\ge 0$, and the first integral
\[
H_2(x, y)=\frac{y^2}{2} +\frac{\e B_1}{4} y^4 + \frac{(1-\e
C_1)x^2}{2}
 \]
if $y\le 0$. Let arbitrary $|\gamma|>0$ and $ R_0>0$  such that
$H_1(\gamma, 0)=R_0$. From this last equality we get $\gamma_{\pm}=\pm
\sqrt{2R_0/(1-\e c_1)}$. Substituting $x=\gamma_{\pm}, y=0$ into $H_2(x,
y)$, we have that $H_2(\gamma_+, 0)=H_2(\gamma_-, 0)$, as shown in
Figure \ref{Center1}. Notice that the origin $O$ is the unique
singularity  of system \eqref{I-piece2} if $|\e|$ is small enough.
Therefore we have a global center at the origin and consequently the
infinity is a center if  $-b_1\e<0$, $-B_1\e<0$ and $ b_1c_1 B_1C_1
\ne0$.

Statement (iii) is  proved and  the proof of  Theorem \ref{th2} is
completed.

\begin{figure}[htb]
\centering
\includegraphics[scale=.30]{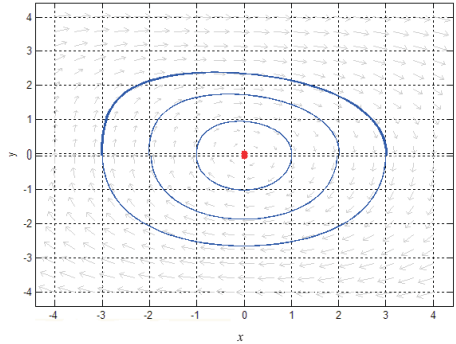}
\caption{   Existence of global center for system \eqref{I-piece2}.  }
\label{up_down1}
\end{figure}

We can illustrate the existence of a global center at the origin in
Theorem \ref{th2} by taking $\e=0.1, b_1=2$, $c_1=0.6$, $d_1=-1$,
$B_1=1$ and $C_1=-0.6$, as it is shown in Figure \ref{up_down1}.

From the proof of  statement (iii) in  Theorem \ref{th2} and the
averaging theory, we have the following results because
$\tilde{f}_{41}(r)$ or  $\tilde{f}_{42}(r)$ has at most two positive
zeros.

\begin{proposition}\label{Pro1}
At most $2$ limit cycles of system \eqref{VII-piece} can bifurcate from the
periodic orbits of the linear system using the averaging theory of
order four, and there are systems \eqref{VII-piece}  with $2$ limit
cycles.
\end{proposition}

\section*{Acknowledgements}

The first author is partially supported by a MINECO grants
MTM2016-77278-P and MTM2013-40998-P, an AGAUR grant number
2014SGR-568, and the grant FP7-PEOPLE-2012-IRSES 318999, and from
the recruitment program of high--end foreign experts of China.

The second author has received funding from the European Union's Horizon 2020 research and innovation
programme under the Marie Sklodowska-Curie grant agreement (No. 655212), and is  partially supported by the National Natural
Science Foundation of China (No. 11431008) and the RFDP of Higher Education of China grant (No.20130073110074).

\end{document}